\documentclass[12pt]{amsart}
\usepackage{a4}
\usepackage{amsmath}
\usepackage{amssymb}

\newtheorem{prop}{\textbf{Proposition}}
\newtheorem{thm}{\textbf{Theorem}}
\newtheorem{lem}{\textbf{Lemma}}

\newenvironment{prove}{\noindent\textbf{\emph{Proof.}}\,}{\hspace*{\fill}\(\boxtimes\)\newline}

\newcounter{numexam}

\newcounter{numrmrk}
\newenvironment{rmrk}{\refstepcounter{numrmrk}\noindent
\textbf{Remark~\thenumrmrk.}}{\mbox{}\(\triangle\)}

\newcommand{\Kk}{\mathcal{K}}
\newcommand{\Pp}{\mathcal{P}}
\newcommand{\Mm}{\mathcal{M}}
\newcommand{\ff}{\varphi}
\newcommand{\cl}[1]{{\overline{#1}}}
\newcommand{\impl}{{\Rightarrow}}
\newcommand{\dimpl}{{\Leftrightarrow}}
\newcommand{\tends}{\underset{n\to\infty}{\longrightarrow}}

\begin{document}

\title[Continuity]
{On the continuity of the Hutchinson operator}

\author[MB{\&}KL]{Michael F. Barnsley, Krzysztof Le\'{s}niak}



\address{Department of Mathematics, 
Australian National University, 
Canberra, ACT, Australia}
\email{michael.barnsley@maths.anu.edu.au}

\address{Faculty of Mathematics and Computer Science, 
Nicolaus Copernicus University,  
Toru\'{n}, Poland}
\email{much@mat.uni.torun.pl}

\begin{abstract}
We investigate whether the Hutchinson operator 
associated with the iterated function system (IFS) is continuous. 
It clarifies several partial results scattered across recent literature.
While the main example for IFS with strict attractor 
was provided by the family of contractions 
(the so-called hyperbolic system), 
the accent was put on ensuring that 
various contractivity-like conditions are preserved 
when the Hutchinson operator is induced, 
unless very recently it was discovered 
that strict attractors are quite often present 
for a large class of noncontractive maps, 
namely projective maps. 
This sets substantial motivation for the study 
whether in general continuity of functions guarantees 
continuity of the induced Hutchinson operator.
\end{abstract}

\maketitle%

\section{Hyperspaces and multifunctions}


We shall assume throughout the paper that \((X,d)\) stands 
for the complete metric space with metric \(d\). 
The closure of \(B\subset X\) will be denoted by \(\cl{B}\).
The \textbf{distance} from point \(b\in X\) to set \(C\subset X\)
is
\[d(b,C) :=\inf_{c\in C} d(b,c).\]
By \(\varepsilon\)\textbf{-neighbourhood} of \(B\subset X\) we understand
\[N_{\varepsilon} B := \{x\in X: d(x,B)<\varepsilon\}.\]
Note that \(N_{\varepsilon} \{b\}\) is an open 
\(\varepsilon\)-ball at \(b\in X\).
The \textbf{Hausdorff distance} between 
\(B\subset X\) and \(C\subset X\) is given by
\[h(B,C) := 
\inf\{\varepsilon >0 : 
B\subset N_{\varepsilon}(C) \wedge 
C\subset N_{\varepsilon}(B)\}.\]

The family of all nonempty subsets of \(X\) is denoted \(\Pp(X)\), 
and becomes an infinite-valued semimetric space when endowed with \(h\) 
(i.e. \(h(B,C)=\infty\) is allowed 
when at least one of the sets \(B,C\) is unbounded, 
and \(h(B,C)=0\) implies only that \(\cl{B}=\cl{C}\)).
The \textbf{hyperspace} of compacta is a metric space \((\Kk(X),h)\), 
where \(\Kk(X)\) consists of nonempty compact subsets of \(X\).
It is complete, since constructing the hyperspace (of closed sets) 
preserves completness and precompact if the base space \(X\) is so
(e.g. \cite{HandbookMulti,Beer,IllanesNadler}). 

Some handy properties of neighborhoods 
and the Hausdorff distance are collected below.
 
\begin{prop}\label{HausdorffNbds}
For \(\emptyset\neq B,B_j,C,C_j\subset X\), \(j\in J\), 
\(\varepsilon > \eta >0\), hold
\begin{enumerate}
\item[(i)] \(\cl{B} = \bigcap_{\varepsilon>0} N_{\varepsilon}B\),
\item[(ii)] \(N_{\varepsilon}\left(\bigcup_{j\in J} B_j\right) = 
\bigcup_{j\in J} N_{\varepsilon}B_j\),
\item[(iii)] \(N_{\varepsilon} \cl{B} = N_{\varepsilon} B\),
\item[(iv)] \(\cl{N_{\eta} B} \subset N_{\varepsilon} B\),
\item[(v)] \(h(B,C)<\varepsilon \impl B\subset N_{\varepsilon}C\),
\item[(vi)] \(B\subset \cl{N_{\varepsilon}C} \wedge
C\subset \cl{N_{\varepsilon}B} \impl h(B,C)\leq \varepsilon\),
\item[(vii)] \(h\left(\bigcup_{j\in J} B_j, \bigcup_{j\in J} C_j\right)
\leq \sup_{j\in J} h(B_j,C_j)\),
\item[(viii)] \(h(\cl{B},\cl{C})=h(B,C)\).
\end{enumerate}
\end{prop}
The proofs are scattered across various sources, e.g.
\cite{HandbookMulti,Beer,IllanesNadler,Superfractals,Wicks,Hutchinson}.

One striking property of compacta we need later is 
the following well-known but very usefull fact from general topology.

\begin{lem}[\cite{Engelking} Thm 4.3.31]\label{LebesgueNo}
Let \(A\subset X\) be a compact set and \(\mathcal{U}\) its open cover,
i.e. a family of open sets with \(\bigcup\mathcal{U}\supset A\).
Then there exists \(\lambda>0\), called a \textbf{Lebesgue number},
such that the family of balls \(\{N_{\lambda}\{a\} : a\in A\}\) 
refines \(\mathcal{U}\), i.e.
\[\forall_{a\in A}\;\exists_{U_{a}\in\mathcal{U}}\;\; 
N_{\lambda}\{a\}\subset U_{a}.\]
\end{lem}


Any map \(\ff : X\to\Pp(X)\) shall be called 
a \textbf{multifunction}. The \textbf{image} of 
\(\emptyset\neq B\subset X\) under \(\ff\) is given by
\[\ff(B) := \bigcup_{b\in B} \ff(b).\]
It is the image of a set via relation \(\ff\subset X\times X\) and 
it should not be confused with the usual image of map
\[R_{\ff}(B) := \{\ff(b):b\in B\} \subset\Pp(X).\]
Obviously \(\ff(B)=\bigcup R_{\ff}(B)\).
Additionally there holds relation \(\ff(\{b\}) = \ff(b)\) 
between the image and value of multifunction at \(b\in B\), 
which also should not lead to ambiguity. 

If the multifunction \(\ff : X\to\Pp(X)\) has compact values, 
then it is written as \(\ff : X\to\Kk(X)\) to remind it. 
Usual function \(f : X\to X\) is identified with 
the multifunction \(\{f\}: X\to\Kk(X)\), 
\(\{f\}(x) := \{f(x)\}\) for \(x\in X\).

The set-theoretic union of multifunctions 
\(\ff_i : X\to \Pp(X)\), \(i\in I\), is defined as 
\(\bigcup_{i\in I} \ff_i : X\to\Pp(X)\),
\(\left(\bigcup_{i\in I} \ff_i\right)(x) := 
\bigcup_{i\in I} \ff_i(x)\) for \(x\in X\). 
The closure of multifunction \(\ff : X\to\Pp(X)\) 
is defined as \(\cl{\ff} : X\to\Pp(X)\),
\(\cl{\ff}(x) := \cl{\ff(x)}\) for \(x\in X\).

We can speak about (uniform) continuity and contractivity 
of the multifunction \(\ff: X\to\Pp(X)\) if we equip \(\Pp(X)\) 
with the Hausdorff distance \(h\). For the readers convenience 
we recall (after \cite{HandbookMulti,Beer,AubinCellina}) 
these definitions. 

A multifunction \(\ff : X\to\Pp(X)\) is
\begin{enumerate}
\item\label{dfcontraction} \textbf{contraction}, if 
there exists a Lipschitz constant 
\(0\leq L<1\)
s.t.
\[\forall_{x_1,x_2\in X}\;\; 
h(\ff(x_1),\ff(x_2))\leq L\cdot d(x_1,x_2)),\]
\item\label{dfweakcontraction} \textbf{weak contraction}, if
there exists a \textit{comparison function} 
\(\eta :(0,\infty)\to[0,\infty)\)
s.t.
\[\forall_{x_1,x_2\in X}\;\; 
h(\ff(x_1),\ff(x_2))\leq \eta(d(x_1,x_2)),\]
\item\label{dfucontinuous} \textbf{uniformly continuous}, if
\[\forall_{\varepsilon > 0}\; \exists_{\delta > 0}\;
\forall_{x_1,x_2\in X}\;\; d(x_1,x_2) < \delta \impl
h(\ff(x_1),\ff(x_2)) < \varepsilon,\]
\item\label{dfcontinuous} \textbf{continuous} at \(x_0\in X\), if
\[\forall_{\varepsilon > 0}\; \exists_{\delta > 0}\;
\forall_{x\in X}\;\; d(x,x_0) < \delta \impl
h(\ff(x),\ff(x_0)) < \varepsilon,\]
\item\label{dfusc} \textbf{upper semicontinuous} at \(x_0\in X\), if
\[\forall_{\varepsilon > 0}\; \exists_{\delta > 0}\;\;
\ff(N_{\delta}\{x_0\}) \subset N_{\varepsilon} \ff(x_0).\]
\end{enumerate}

A comparison function \(\eta:(0,\infty)\to[0,\infty)\)
should satisfy for \(t, t' > 0\):
\begin{enumerate}
 \item[(\(\eta\)-1)] \(t\leq t'\impl \eta(t)\leq \eta(t')\),
 \item[(\(\eta\)-2)] \(\limsup_{r\to t^{+}} \eta(r)<t\),
 \item[(\(\eta\)-3)] \(\lim_{r\to\infty} \left(r - \eta(r)\right) = \infty\).
\end{enumerate}
Directly from the definition it follows that \(\eta\) 
fulfills also:
\(\eta(t)<t\), \(\limsup_{r\to t} \eta(r)<t\) for \(t>0\) 
(e.g. \cite{AFGL}); see \cite{JachymskiWeakContraction,Mate}
for further discussion of conditions put on comparison functions.

\begin{rmrk}
The above continuities are meant in the Hausdorff sense 
and in general should be distincted 
from continuities in the Vietoris sense.
Multivalued contractions are sometimes called 
Nadler contractions.
\end{rmrk}

We have the following hierarchy:
\begin{center}
\begin{tabular}{|c|}\hline
\ref{dfcontraction} \(\impl\) 
\ref{dfweakcontraction} \(\impl\)
\ref{dfucontinuous} \(\impl\)
\ref{dfcontinuous} \(\impl\) 
\ref{dfusc} \\ \hline
\end{tabular}
\end{center}
For example contraction with Lipschitz constant \(L\) 
is weak contraction under \(\eta(t) := L\cdot t\).

One also should remember that upper semicontinuity is 
a notion different from the one used in real analysis;
hence upper hemicontinuity has been also coined 
for multifunctions but in turn to make things worse
it is often reserved for multifunctions which are
upper semicontinuous w.r.t. the weak topology, so we stack
here with a more wide spread term. Upper semicontinuity 
is essentially a multivalued concept, since for 
\(f: X\to X\) and \(\{f\} : X\to\Pp(X)\) it holds that
\(\{f\}\) is upper semicontinuous if and only if it is
continuous which is equivalent to continuity of \(f\).

We remind that the above continuity conditions 
are preserved under taking the closure 
and finite union of multifunctions 
(which can be for example deduced from 
Proposition~\ref{HausdorffNbds}). It is extended by 
Theorem~\ref{compactunioncontinuous} to infinite families
which are compact.


To be able to deal with compact infinite families of maps as 
generalization of finite ones we recall after 
\cite{LesniakBullPAN} the  notion of the 
\textbf{space of multifunctions}. The function space 
\(\Mm(X,X)\) of multifunctions \(\ff : X\to \Pp(X)\) 
shall be equipped with an infinite-valued semimetric,
namely the \textbf{Chebyshev distance} of uniform convergence
\[h_{\sup}(\ff_1,\ff_2) := 
\sup_{x\in X} h(\ff_1(x),\ff_2(x)) =
\sup_{B\in\Pp(X)} h(\ff_1(B),\ff_2(B))\]
for \(\ff_1 , \ff_2 \in\Mm(X,X)\).
Of course \(h_{\sup}(\ff_1,\ff_2)<\infty\) 
for bounded multifunctions \(\ff_i : X\to\Pp(X)\), 
\(i=1,2\), i.e. \(\ff_i(X)\) -- bounded, and 
\(h_{\sup}(\ff_1,\ff_2) = 0\) implies that 
\(\cl{\ff_1} = \cl{\ff_2}\). 
Hence the subspace of bounded (multi)functions 
(with closed values) constitutes standard metric 
function space.

Whenever we speak about (pre)compactness of 
a family of multifunctions we view it in
\((\Mm(X,X),h_{\sup})\).

\begin{thm}[comp. \cite{LesniakBullPAN} Thm.1, \cite{Kieninger} Prop.3.1.3]\label{compactunioncontinuous}
Let \(\ff_i : X\to\Pp(X)\), \(i\in I\), be a 
precompact family of continuous multifunctions. 
Then its union \(\bigcup_{i\in I} \ff_i : X\to\Pp(X)\)
is a continuous multifunction.
\end{thm}
\begin{prove}
Fix \(x_0\in X\), \(\varepsilon >0\) and find a finite
\({\varepsilon}{\slash}{3}\)-net 
\(\{\ff_{j_1},{\ldots},\ff_{j_m}\}\), i.e.
\begin{equation}\label{fihsupnet}
\forall_{i\in I}\;\exists_{k=k(i)\in\{1,{\ldots},m\}}\;
\forall_{x\in X}\;\; h(\ff_i(x),\ff_{j_k}(x)) < \frac{\varepsilon}{3}.
\end{equation}
Since \(\ff_{j_k} : X\to\Pp(X)\), \(k=1,{\ldots},m\), 
are continuous, there exist \(\delta_k >0\) s.t.
\begin{equation}\label{finetcont}
h(\ff_{j_k}(x),\ff_{j_k}(x_0)) < \frac{\varepsilon}{3}
\end{equation}
for \(d(x,x_0) < \delta_k\).

Put \(\delta := \min \{\delta_k : k=1,{\ldots},m\} >0\).
Then according to (\ref{fihsupnet}) and (\ref{finetcont})
we obtain
\begin{eqnarray*}
h(\ff_i(x),\ff_i(x_0))\leq \\
h(\ff_i(x),\ff_{j_{k(i)}}(x)) + 
h(\ff_{j_{k(i)}}(x),\ff_{j_{k(i)}}(x_0)) +
h(\ff_{j_{k(i)}}(x_0),\ff_i(x_0)) < \varepsilon
\end{eqnarray*}
for \(d(x,x_0)< \delta\), \(i\in I\).
Summing up
\[
h\left(\bigcup_{i\in I} \ff_i(x), \bigcup_{i\in I} \ff_i(x_0)\right)
\leq \sup_{i\in I} h(\ff_i(x),\ff_i(x_0)) \leq \varepsilon.
\]
\end{prove}
Similarly a compact family of uniformly continuous multifunctions
[resp. \(L\)-contractions, \(\eta\)-weak contractions] 
has set-theoretic union again of this type. 
Nevertheless if the Lipschitz constant \(L\) 
or the comparison function \(\eta\) is not common 
for all multifunctions in the family , then
we got a serious obstacle (e.g. \cite{Wicks,Kieninger}).

\section{Iterated function systems}

The system \((X, f_i:i\in I)\), consisting of 
a family of maps \(f_i:X\to X\), 
will be called \textbf{iterated function system}, shortly IFS on \(X\). 
When \(I\) is finite we speak about finite IFS. 
Another interesting instance is the compact family of maps 
(compactness understood in the space of functions with the uniform metric).
One can generalize this notion to \textbf{multivalued IFS} 
(e.g. \cite{PetruselRus,Kieninger,LesniakForum,AndresFiser,Ok,AFGL}) 
substituting family of maps by multifunction \(\ff :X\to \Pp(X)\). 
Then any IFS becomes multivalued IFS, if we define 
\(\ff(x) := \{f_i(x) : i\in I\}\), \(x\in X\). 
Of course this leads to abstract investigations 
(which appear for the first time perhaps in \cite{Strother}),
but the benefit of this approach is clarity, scrutiny and unification 
of various types of IFSs, e.g. IFSs with condensation 
and associated inhomogeneous fractals 
from \cite{BarnsleyDemko,Hata} can be cast in this framework 
just by adding one constant multifunction.
Nevertheless the symbolic dynamics cannot be directly used 
in the study of multivalued IFS, unlike for usual IFS, 
since the multifunctions rarely possess decompositions 
into selectors with enough good properties 
(cf. \cite{AubinCellina,HandbookMulti}). 
This is the main drawback in analyzing the structure 
of fractals generated by multifunctions.

The \textbf{Hutchinson operator} \(F:\Pp(X)\to\Pp(X)\), 
associated with a system given by multifunction 
\(\ff: X\to \Pp(X)\) is defined as 
\[F(B) := \cl{\bigcup_{b\in B} \ff(b)}\]
for \(B\in\Pp(X)\). In the case of IFS \((X, f_i:i\in I)\)
this means that
\[F(B) = \cl{\bigcup_{i\in I} f_i(B)} = 
\cl{\{ f_i(b) : i\in I, b\in B \}}\]
for \(B\in\Pp(X)\).

The \(n\)-fold composition of \(F\) is written as \(F^{n}\).
Whenever we speak about abstract IFS, possibly multivalued, the letter
\(F\) denotes its associated Hutchinson operator.

The most important instance of the Hutchinson operator 
is its restriction to the hyperspace of compacta 
\(F:\Kk(X)\to\Kk(X)\), since usually the hyperspace \(\Kk(X)\) 
is perceived as habitat for fractals generated by IFSs.
To be more precise one has to assume that \(F\) sends 
compacta onto compacta. 
Indeed this is fulfilled, when the system 
\(\{f_i\}_{i\in I}\) consists of continuous maps and is finite,
or more generally compact (in the space of functions equipped 
with the metric of uniform convergence).
Still more general condition can be provided for multivalued IFSs.

\begin{prop} 
Let \(\ff : X\to\Kk(X)\) be an upper semicontinuous multifunction 
with compact values.
Then the induced Hutchinson operator \(F : \Pp(X)\to \Pp(X)\) 
transforms compacta into compacta. 
In particular the restriction \(F: \Kk(X)\to\Kk(X)\) 
is well-defined and
\[F(B) = \ff(B)\]
for \(B\in\Kk(X)\).
\end{prop}
\begin{prove}
It is well known that under our assumptions
the image of a compact set is again compact (\cite{HandbookMulti,Beer}).
\end{prove}

Observe yet that \(\ff\) and its closure \(\cl{\ff}\) yield 
the same operator Hutchinson operator 
(\(\cl{\bigcup_{b\in B} \ff(b)} = \cl{\bigcup_{b\in B} \cl{\ff(b)}}\)
for \(B\in\Pp(X)\)), so one can always assume that the values 
of considered multifunctions are closed sets.

\section{Attractors and continuity}

Let \(F:\Pp(X)\to\Pp(X)\) be the Hutchinson operator 
induced by (possibly multivalued) IFS.
We shall say that a compact nonempty set \(A\subset X\) is 
\begin{enumerate}
\item a \textbf{strict attractor}, following 
\cite{AndresFiser,ProjectiveIFS}, when there exists 
an open neighbourhood \(U(A)\supset A\) 
(called basin of attraction) s.t.
\[F^n(B)\tends A,\]
w.r.t. \(h\) for all nonempty compact \(U(A)\supset B\in\Kk(X)\),
\item a \textbf{global maximal attractor}, following
\cite{McGehee,Akin,MelnikValero,Caraballo,LesniakSlovaca},
when 
\begin{equation}\label{DfMaxAttractor}
\forall_{\varepsilon > 0}\;\exists_{n_0}\;
\forall_{n\geq n_0}\;\; F^n(B)\subset N_{\varepsilon} A
\end{equation}
for all nonempty \(B\subset X\), and the set \(A\) 
is minimal w.r.t. the property (\ref{DfMaxAttractor}).
\end{enumerate}
\begin{rmrk}
In a compact space \(X\) strict attractor \(A\) 
with full basin \(U(A)=X\) is also global maximal attractor. 
\end{rmrk}

Only very recently it has been observed that 
aside standard hyperbolic (= contractive) and weakly
contractive IFSs (\cite{Hutchinson,Hata,AndresFiser})
also elliptic (= projective) IFSs possess 
unique strict attractor; see \cite{ProjectiveIFS}.

It can be proved that a global maximal attractor of IFS 
given by upper semicontinuous multifunction is invariant 
(Proposition~5 \cite{LesniakSlovaca}; 
see it also for the case of noncompact closed attractor). 
In the same vein

\begin{prop}
A strict attractor \(A\) of the IFS given by 
an upper semicontinuous multifunction \(\ff : X\to\Pp(X)\)
is invariant i.e. \(F(A)=A\), where 
\(F:\Pp(X)\to\Pp(X)\) is the induced Hutchinson operator.
\end{prop}
\begin{prove}
Fix \(\varepsilon > 0\). By Proposition~2 in 
\cite{LesniakSlovaca} (comp. Lemma~\ref{Cantor2} further)
we know that for some \(\delta >0\)
\[\ff(N_{\delta} A) \subset N_{\frac{\varepsilon}{2}} \ff(A),\]
so
\begin{equation}\label{Fepsde}
F(N_{\delta} A) \subset N_{\varepsilon} F(A).
\end{equation}
From the definition of strict attractor there exists \(n_0\) 
s.t. \(h(F^n(A),A) < \delta\) for \(n\geq n_0\),
so
\begin{equation}\label{Fattract}
F^n(A) \subset N_{\delta} A.
\end{equation}
Combining (\ref{Fepsde}) and (\ref{Fattract}) gives for \(n\geq n_0\)
\[F^{n+1}(A)\subset F(N_{\delta} A) 
\subset N_{\varepsilon} F(A).\]
Thus
\[A=\lim_{n\to\infty} F^{n+1}(A) \subset 
N_{2\varepsilon} F(A)\]
and since \(\varepsilon\) was arbitrary \(A\subset F(A)\).
Due to monotonicity of \(F\) then
\[A\subset F(A)\subset {\ldots} \subset F^n(A) \tends A 
= \cl{\bigcup_{n=1}^{\infty} A_n},\]
which means \(F(A)=A\).
\end{prove}

We want to emphasize that 
if the Hutchinson operator \(F\) is continuous, 
then the invariance of the strict attractor is immediate:
\[A = \lim_{n\to\infty} F^{n+1}(A) = 
F\left(\lim_{n\to\infty} F^n(A)\right) = F(A).\]
This great simplification may be added to a list of 
typical general interest arguments 
supporting the search for conditions under which 
\(F\) is continuous.
It is known that even very simple 
upper semicontinuous multifunctions on compact spaces need
not induce continuous Hutchinson operator 
(\cite{AndresFiser} Counter-Example 1, \cite{Kieninger} Prop.1.5.3).
Another reason to establish continuity of \(F\) 
provides \cite{ChaosGame}.

\section{Known results}

This review section is based on the carefull study of 
articles \cite{Mate,Kieninger,LesniakSlovaca,AndresFiser,AFGL} 
supported by Theorem~\ref{compactunioncontinuous} 
and Proposition~\ref{necessarycondition}.

Let \(\ff : X\to\Pp(X)\) be a multifunction and 
\(F : \Pp(X)\to\Pp(X)\) its associated Hutchinson operator. 
We gather below informations how the continuity of \(\ff\) is preserved 
when inducing \(F\).

\begin{center}
\begin{tabular}{|r|c|}\hline
\mbox{ } &
\(\ff \stackrel{\mbox{?}}{\longleftrightarrow} F\) \\ \hline
contraction & \(\dimpl\) \\
weak contraction & \(\dimpl\) \\
uniformly continuous & \(\dimpl\) \\
continuous & \(\Leftarrow\) \\ \hline 
\end{tabular}
\end{center}
Additionally upper semicontinuity of a multifunction 
\(\ff :X\to \Kk(X)\) with (pre)compact values is equivalent
to upper Vietoris continuity of \(F\).

Compact infinite system of (multi)functions, 
due to Theorem~\ref{compactunioncontinuous} 
with accompanying comments, 
can be turned into a system generated by a single multifunction 
and various continuity conditions are preserved during this process 
as shown in the following table:

\begin{center}
\begin{tabular}{|r|c|}\hline
\mbox{ } &
\(\ff_i (i\in I) \stackrel{\mbox{?}}{\longrightarrow} 
\bigcup_{i\in I} \ff_i\) \\ \hline
\(L\)-contraction  & \(\impl\) \\
\(\eta\)-weak contraction & \(\impl\) \\
uniformly continuous & \(\impl\) \\
continuous & \(\impl\) \\ \hline 
\end{tabular}
\end{center}
The family of multifunctions \(\ff_i : X\to\Pp(X)\), \(i\in I\),
is assumed to be (pre)compact.

In principle both tables cover all situations met 
in applications until very recently (\cite{ProjectiveIFS}). 
The problems with continuity apparent from the first table 
are solved in the main section of this article.

Finally we should also note that the Hutchinson operator 
is often \textit{order-continuous} w.r.t. the inclusion \(\supset\) 
(see \cite{Hayashi,JachymskiGajekPokarowski} 
for details) but from the point of view of applications 
to invariant sets of IFSs \textit{monotonicity} 
together with some amount of compactness suffices 
(which is met in a large class of systems 
as showed in \cite{LesniakForum,Ok,LesniakBullPAN,AFGL}).


\section{Main results}

\begin{prop}[Necessary condition]\label{necessarycondition}
Let \(\ff : (X,d)\to (\Pp(X),h)\) be a multifunction and 
\(F: (\Pp(X),h)\to(\Pp(X),h)\) its induced Hutchinson operator.
If \(F\) is continuous (at singletons), then \(\ff\) is continuous too. 
\end{prop}
\begin{prove}
It is enough to observe that 
\(h(\ff(x),\ff(x_0)) = h(F(\{x\}),F(\{x_0\}))\) for \(x,x_0\in X\).
\end{prove}

\begin{lem}[Cantor--Weierstrass uniform continuity]\label{Cantor}
Let \(\ff : X\to \Pp(X)\) be a continuous multifunction and \(C\subset X\) 
a compact set. Then for every \(\eta>0\)
\[\exists_{\lambda>0}\;\forall_{c\in C}\;
\forall_{z\in N_{\lambda}\{c\}}\;\;
h(\ff(z),\ff(c))<\eta.\] 
\end{lem}
\begin{prove}
For \(\eta{\slash}2 >0\) by continuity of \(\ff\)
with every \(x\in C\) we can associate \(\lambda(x) >0\)
s.t. 
\[\forall_{z\in N_{\lambda(x)}\{x\}}\;
h(\ff(z),\ff(x))<\frac{\eta}{2}.\] 
The open cover 
\(\mathcal{U} := \{ N_{\lambda(c)}\{c\} : c\in C \}\)
of compact \(C\) admits by Lemma~\ref{LebesgueNo} 
a Lebesgue number \(\lambda>0\) i.e.
\[\forall_{c\in C}\; \exists_{x\in C}\;\;
c\in N_{\lambda}\{c\}\subset N_{\lambda(x)}\{x\}.\]
Thus for all \(c\in C\) and \(z\in N_{\lambda}\{c\}\)
\[h(\ff(z),\ff(c))\leq 
h(\ff(z),\ff(x)) + h(\ff(x),\ff(c)) <
\frac{\eta}{2}+\frac{\eta}{2} =\eta.\]
\end{prove}
\begin{rmrk}
The proof can be also performed exactly as the one  
for Theorem 4.3.32 in \cite{Engelking} if we take 
\(\mathcal{U} := \{x\in X : h(\ff(x),\ff(c)) < {\eta}{\slash}{2} \}_{c\in C}\)
to be an open covering of compact \(C\) and \(\lambda >0\) its Lebesgue number.
\end{rmrk}

\begin{lem}[Cantor--Weierstrass uniform continuity II]\label{Cantor2} 
Let \(\ff : X \to \Pp(X)\) be a continuous multifunction and 
\(\emptyset\neq C\subset X\) a nonempty compact set.
Then for every \(\varepsilon >0\) 
there exists \(\delta>0\) such that 
for all \(\emptyset\neq B\subset N_{\delta}\,C\) holds
\[\ff(N_{\delta} B)\subset N_{\varepsilon} \ff(B).\]
\end{lem}
\begin{prove}
Put \(\eta := {\varepsilon}{\slash}{2}\) in Lemma~\ref{Cantor} 
and find \(\lambda >0\) so that
\begin{equation}\label{pointwiseCantor}
h(\ff(z),\ff(c))<\eta
\end{equation}
for all \(c\in C\), \(z\in N_{\lambda}\{c\}\). 
Next define \(\delta := {\lambda}{\slash}{3}\)
and fix nonempty \(B\subset N_{\delta} C\).
We shall verify that for 
\(b\in B\) holds 
\(\ff(N_{\delta} \{b\})\subset N_{\varepsilon} \ff(b)\), 
which in turn gives
\[\ff(N_{\delta}B) = \bigcup_{b\in B} \ff(N_{\delta}\{b\})
\subset \bigcup_{b\in B} N_{\varepsilon}\ff(b) = N_{\varepsilon}\ff(B).\] 

Take \(z\in \ff(N_{\delta} B)\). 
Hence \(z\in \ff(u)\), \(d(u,b)<\delta\)
for some \(b\in B\) and \(u\in N_{\delta} B\). 
Moreover there exists \(c\in C\) satisfying:
\(d(b,c)<\delta <\lambda\) (because \(B\subset N_{\delta} C\))
and \(d(u,c)\leq d(u,b)+d(b,c)< 2\delta < \lambda\). 
Applying (\ref{pointwiseCantor}) yields
\begin{equation*}
h(\ff(u),\ff(b))\leq 
h(\ff(u),\ff(c))+h(\ff(c),\ff(b))< 2\eta = \varepsilon,
\end{equation*}
so \(z\in \ff(u)\subset N_{\varepsilon}\ff(b)\).
\end{prove}

\begin{thm}
Let \(\ff :X\to \Pp(X)\) be a continuous multifunction. 
Then the induced Hutchinson operator \(F:\Pp(X)\to \Pp(X)\) 
is continuous in the Hausdorff metric at every point 
\(C\in \Kk(X)\) in the hyperspace of compacta.
\end{thm}
\begin{prove}
Fix \(C\in \Kk(X)\), \(\varepsilon >0\). 
Find via Lemma~\ref{Cantor2} \(\delta >0\)
such that
\[\ff(N_{\delta} B) \subset N_{\frac{\varepsilon}{2}} \ff(B)\]
for \(B\subset N_{\delta} C\), \(B\in \Pp(X)\).
Thus
\begin{equation}\label{setCantor}
F(N_{\delta} B) =
\cl{\ff(N_{\delta} B)} \subset  
\cl{N_{\frac{\varepsilon}{2}} \ff(B)} \subset 
N_{\varepsilon} \ff(B) \subset 
N_{\varepsilon} F(B)
\end{equation}
for all \(B\subset N_{\delta} C\). In particular,
for every nonempty \(B\subset X\), if \(h(B,C)<\delta\),
then \(C\subset N_{\delta} B\), \(B\subset N_{\delta} C\).
Hence due to (\ref{setCantor}) 
\(F(C)\subset F(N_{\delta} B)\subset N_{\varepsilon} F(B)\),
\(F(B)\subset F(N_{\delta} C)\subset N_{\varepsilon} F(C)\).
Altogether \(h(F(B),F(C))\leq \varepsilon\).
\end{prove}

From the above, due to Theorem~\ref{compactunioncontinuous} 
we obtain

\begin{thm}
Let \((X, f_i : i\in I)\) be a compact (in particular finite) 
IFS consisting of continuous functions. 
Then the induced Hutchinson operator \(F:\Kk(X)\to \Kk(X)\) 
is continuous.
\end{thm}

\end{document}